\documentclass[11pt]{article}
\usepackage{fullpage}
\usepackage{amsmath, amssymb, graphicx, multicol}
\usepackage{xcolor}

\begin{document}

\begin{center}
\textbf{APPLICATIONS OF TEACHING SECONDARY MATHEMATICS IN UNDERGRADUATE MATHEMATICS COURSES} 
\vskip0.25in
Elizabeth G. Arnold, Elizabeth A. Burroughs, Elizabeth W. Fulton, James A. Mendoza \'Alvarez \\
Colorado State University, Montana State University, Montana State University, University of Texas at Arlington \\
\vskip0.25in
\color{blue}{Paper for TSG33 at the 14th International Congress on Mathematical Education (2020). Shanghai, China. [Delayed to July 2021 due to COVID-19 pandemic]}
\end{center}

\noindent\textit{Robust preparation of future secondary mathematics teachers requires attention to the acquisition of mathematical knowledge for teaching. Many future teachers learn mathematics content primarily through mathematics major courses that are taught by mathematicians who do not specialize in teacher preparation. How can mathematics education researchers assist mathematicians in making explicit connections between the content of undergraduate mathematics courses and the content of secondary mathematics? We present an articulation of five types of connections that can be used in secondary mathematics teacher preparation and give examples of question prompts that mathematicians can use as applications of teaching secondary mathematics in undergraduate mathematics courses.
}

\section*{INTRODUCTION}
Secondary mathematics teacher preparation routinely includes the study of both mathematics content and pedagogical strategies. It is well understood that secondary mathematics teachers must have a strong understanding of the mathematics taught in secondary schools (Ferrini-Mundy \& Findell, 2001). The Mathematical Education of Teachers II (MET II) Report (CBMS, 2012) recommends that prospective secondary teachers complete the equivalent of an undergraduate major in mathematics with a focus on examining secondary school mathematics from an advanced perspective, emphasizing the importance of making connections between the mathematics undergraduates are learning and the school mathematics they will be teaching. In this vein, Wasserman (2018b) recommends that the content of secondary mathematics ought to inform how undergraduate mathematics courses are taught. Despite these recommendations, many undergraduate mathematics major programs in the United States do not meet these needs: ``Evidence suggests that prospective high school mathematics teachers, who earn a mathematics major or its equivalent, do not have sufficiently deep understanding of the mathematics of the high school curriculum" (Speer, King \& Howell, 2015, p. 107). Zazkis and Leikin (2010) found that ``many teachers perceive their undergraduate studies of mathematics as having little relevance to their teaching practice" (p. 1). There is a gap between the mathematics that prospective teachers are learning and the school mathematics they will teach, pointing to the need to focus on connections to school mathematics within the undergraduate mathematics curriculum. 

Research in mathematics education has shown that it is not enough for future teachers to be proficient with school mathematics; future teachers must also develop a kind of specialized knowledge for teaching. Mathematical knowledge for teaching (MKT) is a theory developed by Ball, Thames, and Phelps (2008), focused on elementary mathematics teaching, that encompasses the knowledge teachers require in order to ``carry out the work of teaching mathematics" (p. 395) and plays an important role in teacher preparation programs. Stylianides and Stylianides (2014) articulate the need for teacher preparation programs to incorporate specific knowledge needed for teaching as well as opportunities for prospective teachers to apply this knowledge in context of their future practice.  Yet, in the United States, it is often the case that prospective secondary mathematics teachers enroll in traditional mathematics major courses that are often taught by mathematicians who do not specialize in teacher preparation nor have a robust understanding of mathematical knowledge for teaching.

This research operationalizes mathematical knowledge for teaching to support mathematicians in addressing it in their instruction. We present practical implications for prospective mathematics teacher preparation that result from an articulation of five types of connections that can support all instructors of mathematics content courses.

\section*{THEORETICAL CONSIDERATIONS}
\begin{figure}[ht]
\centering
\includegraphics{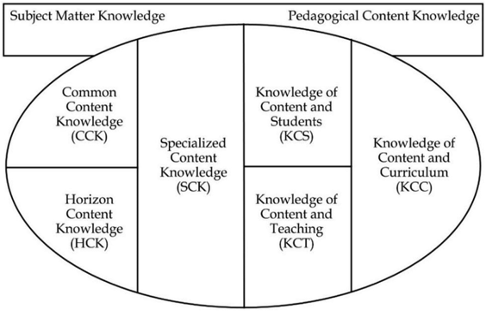}
\caption{ \label{fig1}The components of mathematical knowledge for teaching (Ball et al., 2008)}
\end{figure}

Teaching is a complex and multifaceted skill that employs a wide range of knowledge. Shulman (1986) was among the first to propose that in order to teach mathematics effectively, teachers must possess more than a separate understanding of content and pedagogy. He described another component of teaching, called pedagogical content knowledge, in which teachers combine their pedagogical knowledge with their content knowledge. This idea has been influential in teacher education research, and Ball et al. (2008) extended the idea of pedagogical content knowledge when developing the construct of mathematical knowledge for teaching (Figure~\ref{fig1}).

Ball et al. (2008) address two broad domains in the construct of MKT: subject matter knowledge and pedagogical content knowledge. Subject matter knowledge comprises a teacher’s knowledge of content within a course, how this content connects to previous and future topics, and different ways content can be interpreted to assist in teaching. Pedagogical content knowledge comprises a teacher’s knowledge of the content in relation to teaching practices, student learning, and the curriculum. Six kinds of knowledge relevant to teaching elementary mathematics make up these two broad domains.

\begin{itemize}
	\item \textit{Common Content Knowledge} is the mathematical knowledge and skills used in many situations (teaching and non-teaching). It includes, among other things, being able to solve mathematical problems, recognize errors, and use terms and notation correctly. 
	\item \textit{Horizon Content Knowledge} is knowledge of the mathematics that follows or could follow the mathematics being taught. It is an awareness of how the mathematical topics span a curriculum. 
	\item \textit{Specialized Content Knowledge} is the mathematical knowledge and skills (beyond common content knowledge) that are specific to teaching. 
	\item \textit{Knowledge of Content and Students} is a knowledge about how students learn mathematics and being able to anticipate common student errors.
	\item \textit{Knowledge of Content and Teaching} is knowledge about mathematical tasks and teaching, such as sequencing and selecting examples to help students learn.
	\item \textit{Knowledge of Content and Curriculum} is the knowledge of curriculum for a given topic. This includes being familiar with materials and resources for teaching mathematics. 
\end{itemize}

Other researchers (e.g., Wasserman, 2018a; Wasserman et al., 2019) have addressed ways to develop mathematical knowledge for teaching in the preparation of secondary mathematics teachers. Wasserman et al. (2019) articulate four ways in which prospective teachers can engage in developing mathematical knowledge for teaching and describe these as existing along a spectrum spanning from ``mathematical in nature" to ``pedagogical in nature" (p. 821). \textit{Content} tasks focus on the undergraduate mathematics that serves to deepen prospective teachers’ understanding of secondary mathematics. \textit{Modeled Instruction} consists of different modes of instruction implemented by an instructor (e.g., inquiry-based learning) that can influence how prospective teachers will teach. Mathematical practices and mathematical habits of mind are the foundation for \textit{Disciplinary Practice}; Wasserman (2018a) describes this as focusing on what it means to ``do" mathematics. Lastly, \textit{Classroom Teaching} tasks apply content to a specific teaching situation such as designing problems, sequencing activities, and describing, for example, what it means to multiply two fractions.

\section*{FIVE CONNECTIONS TO TEACHING SECONDARY MATHEMATICS}
Using Ball et. al.’s (2008) delineation of the six categories of MKT as our theoretical base, we defined five types of ``connections for teaching" between undergraduate-level mathematics and knowledge for teaching secondary mathematics (Table~\ref{tab1}). The MET II Report (CBMS, 2012) suggests engaging prospective teachers in understanding ``the connections" (p. 54) between undergraduate mathematics and school mathematics, but leaves these connections undefined. Though education researchers generally understand that the knowledge needed for teaching encompasses more than content knowledge, there is not such awareness among non-specialists, that is, among those who teach the bulk of the undergraduate courses taken by prospective teachers but are not educational researchers, leaving a compelling need for an articulation of how to make such connections. 

The connections we formulate here complement the existing literature and extend the construct of MKT in a way that is useful for mathematicians who have a role in preparing secondary teachers. When developing our five connections, we first considered each broad domain of MKT separately. In considering how to translate subject matter knowledge for teaching secondary mathematics, we asked ourselves the following question: what should undergraduate prospective teachers understand about the undergraduate mathematics that is essential to the study of school mathematics? We determined that understanding mathematics at this level requires being able to: 1) solve mathematical problems foundational to school mathematics; 2) explain the reasoning supporting the mathematical concepts; and 3) know what mathematics comes before and after the topic at hand. We used these three ideas to formulate our first three connections. In considering what it means for an undergraduate to develop pedagogical content knowledge in a mathematics course, we asked ourselves: what do secondary students understand, and how can undergraduates practice helping these students with their understanding? We arrived at the importance of undergraduate prospective teachers’ ability to: 4) look at student work and recognize and articulate what a student does and does not understand mathematically; and then 5) think of questions that will help guide that student’s understanding. These ideas comprise our last two connections.

\begin{table}[ht]
  \begin{center}
    \caption{Five connections to teaching secondary mathematics.}
    \label{tab1}
    \begin{tabular}{|l|l|} \hline
    \textbf{Connection} &  \textbf{Description} \\ \hline
    1. \textit{Content Knowledge} & Undergraduates use course content in applications \\
    & or to answer mathematical questions in the course. \\ \hline
    2. \textit{Explaining Mathematical Content} & Undergraduates justify mathematical procedures \\
    & or theorems and use of related mathematical \\
    & concepts. \\ \hline
    3. \textit{Looking Back / Looking Forward} & Undergraduates explain how mathematics topics \\
    & are related over a span of K-12 curriculum \\
    & through undergraduate mathematics.\\ \hline
    4. \textit{School Student Thinking} & Undergraduates evaluate the mathematics \\
    & underlying a student’s work and explain what\\
    &  that student may understand. \\ \hline
    5. \textit{Guiding School Students’ Understanding} & Undergraduates pose or evaluate guiding questions \\
    & to help a hypothetical student understand a \\
    & mathematical concept and explain how the \\
    & questions may guide the student’s learning. \\ \hline
        \end{tabular}
  \end{center}
\end{table}

We note that we do not intend for these five connections to be a direct mapping to the kinds of knowledge described in the MKT framework. Rather, we formulate these connections in this way to facilitate addressing them in curricular materials used by mathematicians when preparing future secondary teachers. Additionally, we do not claim that these five connections are disjoint; we fully expect that some problems and activities will engage undergraduates in more than one of these connections.

We illustrate these five connections with examples. We created instructional materials that incorporate these five types of connections into lessons on a selection of topics in undergraduate mathematics courses, and we have field tested these instructional materials in a variety of settings. We selected lessons in the undergraduate content areas of abstract algebra, calculus, discrete mathematics, and statistics, because these areas are foundational for the mathematics studied in secondary schools.

\subsection*{Content Knowledge}
\textit{Content Knowledge} connections require undergraduates to use content they learn in a course to answer mathematical questions, where the undergraduate content directly relates to the content of secondary school mathematics. The following is a problem from a Statistics lesson on margin of error:

\begin{figure}[ht]
\centering
\includegraphics[scale=0.85]{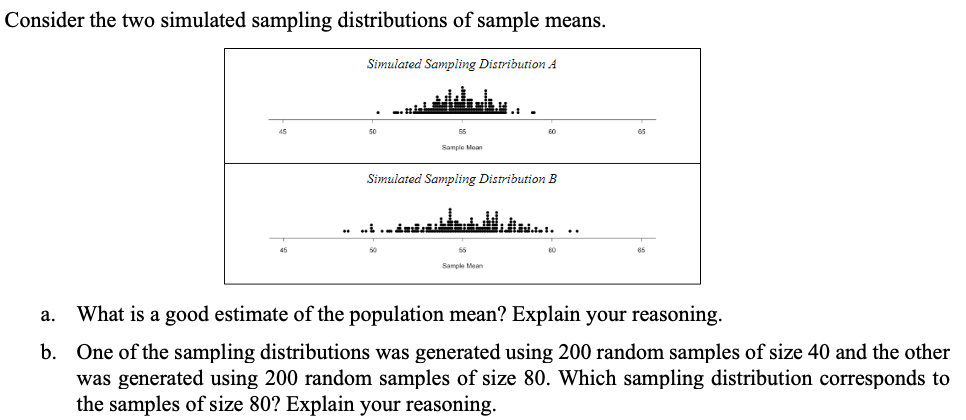}
\end{figure}

This problem engages undergraduates in using data from a random sample to estimate a population mean with the primary purpose of enhancing their understanding of the relationship between sample size and margin of error, concepts they will teach their future students. These types of problems are similar to those most mathematicians normally use in teaching undergraduate courses; we identify them as Content Knowledge Connections because we have specifically chosen problems that focus on concepts that are firmly grounded in secondary content.

\subsection*{Explaining Mathematical Content}
\textit{Explaining Mathematical Content} connections require undergraduates to justify a mathematical claim, procedure, formula or theorem. Going beyond content knowledge, this type of connection enforces deep understanding of content. An example from a lesson on groups of transformations for an Abstract Algebra course is:

\begin{figure}[ht]
\centering
\includegraphics[scale=0.85]{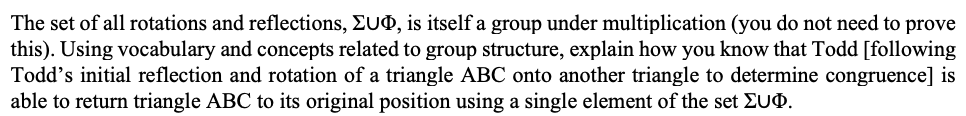}
\end{figure}

This problem focuses on undergraduates’ understanding of groups and the logical implications that follow from a given algebraic structure. These types of problems involve important communication skills which are particularly valuable for prospective teachers whose jobs will require them to explain mathematical concepts to others. The lessons we developed provide problems and questions as a resource for mathematicians who may not normally require undergraduates to reflect and discuss why a mathematical approach is valid.

\subsection*{Looking Back/Looking Forward}
\textit{Looking Back/Looking Forward} connections require undergraduates to examine how the mathematics topics they are studying relate to topics within the secondary curriculum. The following is a problem from a Calculus lesson on derivatives of inverse functions:

\begin{figure}[ht]
\centering
\includegraphics[scale=0.85]{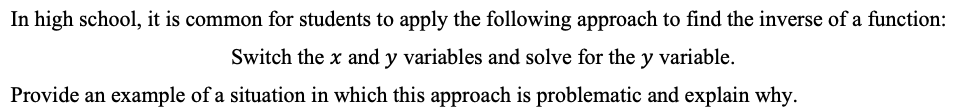}
\end{figure}

These types of questions situate mathematics topics in the context of secondary teaching and learning. This connection differs from \textit{Content Knowledge} connections because, when developing \textit{Content Knowledge} problems, the relationship between the college mathematics to school mathematics is considered ``behind the scenes"; the instructor identifies the mathematics in undergraduate courses that is important for teaching secondary mathematics and creates problems that specifically address the mathematics without explicitly referencing the connection. When writing \textit{Looking Back/Looking Forward} problems, the relationship between the college and secondary mathematics is made explicit through the context in which the problem is posed. The lessons we developed support mathematicians who may not be familiar enough with the secondary curriculum to include such \textit{Looking Back/Looking Forward} connections on their own.

\subsection*{School Student Thinking}
\textit{School Student Thinking} connections require undergraduates to evaluate the mathematics presented in a hypothetical student’s work and explain what that student may understand. The following is a problem from a lesson on the Foundations of Divisibility in an introduction to proofs course:

\begin{figure}[ht]
\centering
\includegraphics[scale=0.85]{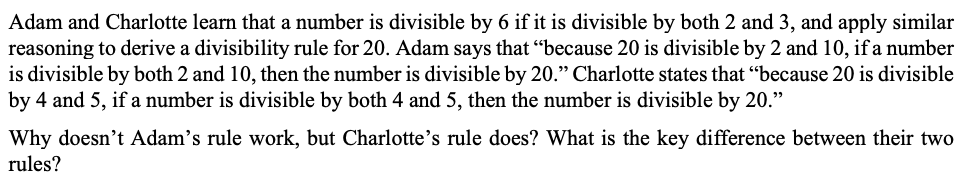}
\end{figure}

Questions like these give prospective teachers opportunities to recognize mathematical understanding in incomplete or incorrect solutions. Questions addressing \textit{School Student Thinking} also prompt undergraduates to explain why students may use different, yet correct, approaches to a problem. The questions can use student thinking to guide undergraduates through unknown mathematics. These types of questions provide opportunities for undergraduates to use or develop their mathematical understanding by investigating other students’ perspectives. Our lessons provide explicit curriculum support for mathematicians by providing examples of common school student errors with prompts for undergraduates to consider.

\subsection*{Guiding School Students’ Understanding} 
\textit{Guiding School Students’ Understanding} connections require undergraduates to pose or evaluate questions to help a hypothetical student understand a mathematical concept and to explain how the posed questions may guide student learning. The following is a problem from a Binomial Theorem lesson in an introduction to proofs course:

\begin{figure}[ht]
\centering
\includegraphics[scale=0.85]{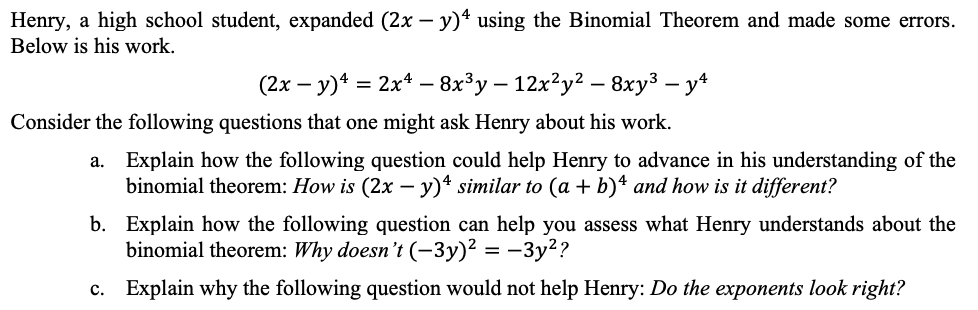}
\end{figure}

Posing mathematical questions that probe understanding or that advance understanding is an important skill for future teachers to learn and practice (Smith and Stein, 2018). Mathematicians who teach undergraduates may not be familiar enough with the school curriculum to provide such opportunities without explicit curricular support. 

\section*{DISCUSSION}
Robust preparation of future secondary mathematics teachers requires attention to the acquisition of mathematical knowledge for teaching. We offer these five connections as a means to tailor the content in undergraduate mathematics courses to the work of teaching. Our particular interest is how to incorporate connections into mathematics courses commonly included in a mathematics major for prospective secondary mathematics teachers taught by non-specialists in mathematics education research. In turn, this will provide prospective secondary mathematics teachers with multiple opportunities to build aspects of their own mathematical knowledge for teaching throughout their undergraduate mathematics studies. 

The five connections we describe align with Wasserman’s (2018a) delineation of the spectrum of ideas that help make the study of undergraduate mathematics content relevant to future secondary teachers. Our first three connections, situated primarily in content knowledge, align with the kinds of curriculum decisions most mathematics instructors already make. Focusing on the content of mathematics, both \textit{Content Knowledge} and \textit{Explaining Mathematical Content} connections are typically incorporated into most undergraduate mathematics courses. The difference in our context is that the content specifically relates to or arises from high school content in some way. Because many mathematicians consider ways to motivate the content of their course by relating the content to previously studied mathematics concepts, the \textit{Looking Back/Looking Forward} connection may enhance ways that mathematicians typically motivate concepts by providing more support for linking to undergraduates’ high school experiences in mathematics or to foundational high school topics prospective secondary mathematics teachers will have to teach.  

Our articulation of the two pedagogically-focused connections (\textit{Student Thinking} and \textit{Guiding Student Understanding}) offers a novel perspective on teacher preparation in content courses. These connections are complementary to, but distinct from, Wasserman’s (2018a) \textit{Classroom Teaching} activities. They are also distinct from the kinds of tasks that mathematics instructors would usually set for their undergraduate students. It is a teacher’s role to evaluate student work and pose questions to assess and advance their students’ understanding. This is challenging work, and prospective teachers need exposure to these practices in their preparation. Problems incorporating \textit{Student Thinking} and \textit{Guiding Student Understanding} connections differ from traditional problems found in undergraduate mathematics courses. Nonetheless, they provide a novel way to further explore undergraduates’ conceptual understanding via the ways in which undergraduates choose to guide student understanding or pose questions to probe student thinking. These connections, in particular, extend the MET II Report’s (CBMS, 2012) call to make connections between the mathematics prospective teachers are learning and the school mathematics they will be teaching by emphasizing how the mathematics in their undergraduate courses can help them understand their future students’ mathematical work. 

These five connections extend the literature on mathematical knowledge for teaching and offer structure for how we can support prospective teachers in content courses. We recommend the inclusion of these five connections across undergraduate courses as a means for mathematicians to make explicit connections between the content of undergraduate mathematics courses and the content of secondary mathematics. We also envision that the five connections may promote broad use of  embedding applications of teaching secondary school mathematics into mainstream mathematics courses for secondary mathematics teachers. 

\subsection*{Additional Information}
This material is based upon work supported by the National Science Foundation under Award No. DUE-1726624. 

\subsection*{References}
Ball, D. L., Thames, M. H., \& Phelps, G. (2008). Content knowledge for teaching: What makes it special? \textit{Journal of Teacher Education}, 59(5), 389–407. \\

\noindent Conference Board of the Mathematical Sciences (2012). \textit{Mathematical Education of Teachers II Report}. Providence RI and Washington DC: American Mathematical Society and Mathematical Association of America.\\ 

\noindent Ferrini-Mundy, J., \& Findell, B. (2001). The mathematical education of prospective teachers of secondary school mathematics: Old assumptions, new challenges. In T. Rishel (Ed.), \textit{Mathematics and the mathematical sciences in 2010: What should students know?} (pp. 31–41). Washington, DC: Mathematical Association of America.\\

\noindent Speer, N. M., King, K. D., \& Howell, H. (2014). Definitions of mathematical knowledge for teaching: Using these constructs in research on secondary and college mathematics teachers. \textit{Journal of Mathematics Teacher Education}, 18(2), 105–122.\\

\noindent Stylianides, A., \& Stylianides, G. (2014). Viewing “mathematics for teaching” as a form of applied mathematics: Implications for the mathematical preparation of teachers. \textit{Notices of the AMS}, 61(3), 266–276.\\

\noindent Shulman, L. (1986). Knowledge and teaching: Foundations of the new reform. \textit{Harvard Educational Review}, 57, 1–22. \\

\noindent Smith, M. S., \& Stein, M. K. (2018). \textit{5 practices for orchestrating productive mathematics discussions.} Reston, VA: National Council of Teachers of Mathematics; Corwin. \\

\noindent Wasserman, N. H. (2018a). Exploring advanced mathematics courses and content for secondary mathematics teachers. In N. H. Wasserman (Ed.), \textit{Connecting Abstract Algebra to Secondary Mathematics for Secondary Mathematics Teachers} (1–16). Switzerland: Springer. \\

\noindent Wasserman, N. H. (2018b). Knowledge of nonlocal mathematics for teaching. \textit{The Journal of Mathematical Behavior}, 49, 116–128.\\

\noindent Wasserman, N. H., Zazkis, R., Baldinger, E., Marmur, O., \& Murray, E. (2019). Points of connection to secondary teaching in undergraduate mathematics courses. \textit{Proceedings of the 22nd Annual Conference on Research in Undergraduate Mathematics Education.} Oklahoma City, OK.\\

\noindent Zazkis, R., \& Leikin, R. (2010). Advanced mathematical knowledge in teaching practice: Perceptions of secondary mathematics teachers. \textit{Mathematical Thinking and Learning}, 12(4), 263–281. \\

\end{document}